# On the number of collisions in beta(2, *b*)-coalescents

ALEX IKSANOV[1,*], ALEX MARYNYCH[1,**] and MARTIN MÖHLE[2]

[1]*Faculty of Cybernetics, National T. Shevchenko University, 01033 Kiev, Ukraine.*
*E-mails: [*]iksan@unicyb.kiev.ua; [**]marinich@voliacable.com*
[2]*Mathematical Institute, University of Düsseldorf, Universitätsstraße 1, 40225 Düsseldorf, Germany. E-mail: moehle@math.uni-duesseldorf.de*

Expansions are provided for the moments of the number of collisions $X_n$ in the $\beta(2,b)$-coalescent restricted to the set $\{1,\ldots,n\}$. We verify that $X_n/\mathbb{E}X_n$ converges almost surely to one and that $X_n$, properly normalized, weakly converges to the standard normal law. These results complement previously known facts concerning the number of collisions in $\beta(a,b)$-coalescents with $a \in (0,2)$ and $b=1$, and $a>2$ and $b>0$. The case $a=2$ is a kind of 'border situation' which seems not to be amenable to approaches used for $a \neq 2$.

*Keywords:* asymptotics of moments; beta-coalescent; number of collisions; random regenerative composition; recursion with random indices

## 1. Introduction and main results

Let $\mathcal{E}$ denote the set of all equivalence relations (partitions) on $\mathbb{N}$. For $n \in \mathbb{N}$, let $\varrho_n : \mathcal{E} \to \mathcal{E}_n$ denote the natural restriction to the set $\mathcal{E}_n$ of all equivalence relations on $\{1,\ldots,n\}$. For $\xi \in \mathcal{E}_n$ let $|\xi|$ denote the number of blocks (equivalence classes) of $\xi$.

Pitman [15] and Sagitov [17] independently introduced coalescent processes with multiple collisions. These Markovian processes with state space $\mathcal{E}$ are characterized by a finite measure $\Lambda$ on $[0,1]$ and hence are also called $\Lambda$-coalescent processes. For a $\Lambda$-coalescent $\{\Pi_t : t \geq 0\}$, it is known that

$$g_{nk} := \lim_{t \searrow 0} \frac{\mathbb{P}\{|\varrho_n \Pi_t| = k\}}{t} = \binom{n}{k-1} \int_{[0,1]} x^{n-k-1}(1-x)^{k-1} \Lambda(\mathrm{d}x) \qquad (1)$$

for $k, n \in \mathbb{N}$ with $k < n$. Let

$$g_n := \lim_{t \searrow 0} \frac{1 - \mathbb{P}\{|\varrho_n \Pi_t| = n\}}{t} = \sum_{k=1}^{n-1} g_{nk}, \qquad n \in \mathbb{N},$$

denote the total rates. Recently, there appeared several papers [2, 3, 4, 6, 8, 9, 10] dealing with certain functionals of the restricted coalescent process $\{\varrho_n \Pi_t : t \geq 0\}$ (for some particular choices of $\Lambda$). Functionals under consideration in these papers are (i) the number $X_n$ of collision events (jumps) that take place until there is just a single block, and (ii) the total branch length $L_n$, that is, the sum of the length of all branches of the restricted coalescent tree. Such functionals are important for biological and statistical applications because they are closely related to the number of mutations on the restricted coalescent tree, if it is assumed that mutations occur independently of the underlying genealogical tree (neutrality) on each branch of the tree according to some homogeneous Poisson process with parameter $r > 0$ (coalescent with mutation).

In particular, the weak asymptotic behavior of the number of collisions $X_n$ is known for $\beta(a,b)$-coalescents with $a \in (0,2)$ and $b = 1$, and $a > 2$ and $b > 0$. We briefly recall the corresponding weak convergence results because they provide insight into the role of the parameter $a$ of the beta distribution $\Lambda = \beta(a,b)$ in this model.

If $0 < a < 1$ and $b = 1$, then (see [10])

$$\frac{X_n - n(\alpha - 1)}{(\alpha - 1)n^{1/\alpha}} \overset{d}{\to} X,$$

where $\alpha := 2 - a$ and $X$ is an $\alpha$-stable random variable with characteristic function $\mathbb{E} e^{itX} = \exp(|t|^\alpha(\cos(\pi\alpha/2) + i\sin(\pi\alpha/2)\,\text{sgn}(t)))$, $t \in \mathbb{R}$. Gnedin and Yakubovich ([8], Theorem 9) used analytic methods to generalize this result to $\Lambda$-coalescents satisfying $\Lambda([0,x]) = Ax^a + O(x^{a+\zeta})$ as $x \downarrow 0$, where $a \in (0,1)$, $A > 0$ and $\zeta > \max\{(2-a)^2/(5 - 5a + a^2), 1 - a\}$.

If $a = b = 1$ (Bolthausen–Sznitman coalescent), then (see [4, 9, 10])

$$\frac{(\log n)^2}{n} X_n - \log(n \log n) \overset{d}{\to} X,$$

where $X$ is a 1-stable random variable with characteristic function $\mathbb{E} e^{itX} = \exp(it \log|t| - \frac{\pi}{2}|t|)$, $t \in \mathbb{R}$.

If $1 < a < 2$ and $b = 1$, then (see [10])

$$\frac{X_n}{\Gamma(2 - \alpha)n^\alpha} \overset{d}{\to} \int_0^\infty e^{-U_t}\,dt,$$

where $\alpha := 2 - a$ and $\{U_t : t \geq 0\}$ is a drift-free subordinator with Lévy measure $\nu(dt) = e^{-t/\alpha}/((1 - e^{-t/\alpha})^{\alpha+1})\,dt$, $t > 0$.

If $a > 2$ and $b > 0$, then (see [6])

$$\frac{X_n - \mu_1^{-1} \log n}{(\mu_2 \mu_1^{-3} \log n)^{1/2}} \overset{d}{\to} X,$$

where $X$ is a random variable with the standard normal law, $\mu_1 := \Psi(a - 2 + b) - \Psi(b)$, $\mu_2 := \Psi'(b) - \Psi'(a - 2 + b)$ and $\Psi(z) := (d/dz) \log \Gamma(z)$ denotes the logarithmic derivative of the gamma function.

There is also very precise information available concerning the asymptotics of the moments of $X_n$ for $\beta(a,1)$-coalescents with $a \in (0,1]$. For more details, we refer to [10] and [14].

The convergence results above indicate, in particular, that the two special parameter values $a = 1$ and $a = 2$ play a kind of threshold role when studying the limiting behavior of $X_n$. This paper focuses on the asymptotics of $X_n$ for $\beta(a,b)$-coalescents with parameter $a = 2$ (and arbitrary $b > 0$). To the best of our knowledge, no convergence results have yet been provided for these particular beta coalescents.

From the structure of the coalescent process, it follows that $\{X_n : n \in \mathbb{N}\}$ satisfies the recursion

$$X_1 := 0 \quad \text{and} \quad X_n \stackrel{d}{=} X_{n-I_n} + 1, \qquad n \in \{2,3,\ldots\}, \tag{2}$$

where $I_n$ is a random variable independent of $X_2, \ldots, X_{n-1}$ with distribution $\mathbb{P}\{I_n = n-k\} = g_{nk}/g_n$, $k \in \{1,\ldots,n-1\}$. The random variable $n - I_n$ is the (random) state of the process $\{|\varrho_n \Pi_t| : t \geq 0\}$ after its first jump.

As already mentioned above, our aim is to investigate the asymptotic behavior of $X_n$ for $\beta(2,b)$-coalescents with $b > 0$. In this case, $I_n$ has distribution

$$\mathbb{P}\{I_n = k\} = \frac{\Gamma(n-k+b-1)\Gamma(n+1)}{(k+1)\Gamma(n-k)\Gamma(n+b)H(n,b)}, \qquad k \in \{1,\ldots,n-1\}, \tag{3}$$

where

$$H(n,b) := \frac{b}{b+n-1} + \Psi(b+n-1) - \Psi(b) - 1, \qquad n \in \mathbb{N}, b > 0.$$

Note that $\Psi(b+n-1) = \log n + \mathrm{O}(1/n)$, $n \to \infty$ (see (6.3.18) in [1]) and therefore

$$H(n,b) = \log n - \Psi(b) - 1 + \mathrm{O}\left(\frac{1}{n}\right), \qquad n \to \infty. \tag{4}$$

In the proofs, we will need the asymptotics of the total rates

$$g_n = \frac{H(n,b)}{B(2,b)} \sim \frac{\log n}{B(2,b)}, \qquad n \to \infty, \tag{5}$$

where $B(x,y) := \int_0^1 u^{x-1}(1-u)^{y-1}\,\mathrm{d}u$, $x, y > 0$, denotes the beta function. Moreover, we will use the Lévy measure $\mu_b$ on $(0,\infty)$ defined via

$$\mu_b(\mathrm{d}t) := \frac{\mathrm{e}^{-bt}}{1-\mathrm{e}^{-t}}\,\mathrm{d}t, \qquad t > 0, b > 0. \tag{6}$$

Note that $\mu_b$ has moments

$$\begin{aligned}
m_r^{(b)} &:= \int_{(0,\infty)} t^r \mu_b(\mathrm{d}t) = \int_{(0,1)} (-\log(1-x))^r \frac{(1-x)^{b-1}}{x}\,\mathrm{d}x \\
&= \Gamma(r+1)\zeta(r+1,b), \qquad r > 0,
\end{aligned} \tag{7}$$

which follows from a Hurwitz identity (see, for example, (23.2.7) in [1]). Here, $\zeta(z,b) = \sum_{i=0}^{\infty}(i+b)^{-z}$, $\text{Re}(z) > 1$, is the Hurwitz zeta function. Our first result presents the asymptotic expansions of the moments of $X_n$. For convenience, we use the notation $\log^k n := (\log(n))^k$, $k, n \in \mathbb{N}$.

**Theorem 1.1 (Expansion of moments).** *As $n \to \infty$, for $k \in \mathbb{N}$,*

$$\mathbb{E}X_n^k = \frac{1}{(2m_1)^k}\log^{2k} n + \frac{2k((2k+1)m_2 + 6cm_1)}{3(2m_1)^{k+1}}\log^{2k-1} n + \mathrm{O}(\log^{2k-2} n),$$

*where $m_1 := m_1^{(b)} = \zeta(2,b)$ and $m_2 := m_2^{(b)} = 2\zeta(3,b)$ (see (7)), and $c := -\Psi(b) - 1$. In particular, the variance $\mathbb{D}X_n$ has the asymptotic expansion*

$$\mathbb{D}X_n = \frac{m_2}{3m_1^3}\log^3 n + \mathrm{O}(\log^2 n) = \frac{2\zeta(3,b)}{3\zeta^3(2,b)}\log^3 n + \mathrm{O}(\log^2 n).$$

*Remark 1.2.* Let $\{S_t : t \geq 0\}$ be a drift-free subordinator with Lévy measure (6). For $n \in \mathbb{N}$, let $Y_n$ ($Z_n$) be the number of parts (with more than one point) of a regenerative composition arising from throwing $n$ independent (random) points, which are independent of $\{S_t : t \geq 0\}$ and all uniformly distributed on $[0,1]$, on the closed range of the multiplicative subordinator $\{1 - e^{-S_t} : t \geq 0\}$.

According to (19) and (22) in [7], $\mathbb{E}Y_n$ and $\mathbb{E}Y_n^2$ admit almost the same asymptotic expansions as $\mathbb{E}X_n$ and $\mathbb{E}X_n^2$, the only difference being that our $c$ equals $-\Psi(b) - 1$ and their $c$ equals $-\Psi(b)$. According to (19) and Theorem 14 in [7], $\mathbb{E}Z_n$ admits exactly the same asymptotic expansion as $\mathbb{E}X_n$. According to (24) in [7], $\mathbb{D}Y_n$ has the same asymptotic expansion as $\mathbb{D}X_n$. These observations strongly suggest that $X_n$ and $Y_n$ may have a similar asymptotic behavior.

*Remark 1.3.* For $t \geq 0$, let $\{f_i(t) : i \in \mathbb{N}\}$ be the sequence (in some order) of the asymptotic frequencies of the random exchangeable partition $\Pi_t$. Note that $\int_{[0,1]} x^{-1}\Lambda(\mathrm{d}x) < \infty$ for $\Lambda = \beta(2,b)$, $b > 0$. Therefore, by Proposition 26 in [15], $\{\widehat{S}_t := -\log(1 - \sum_{i=1}^{\infty} f_i(t)) : t \geq 0\}$ is a version of $\{S_t : t \geq 0\}$. We will come back to this remark later in the proofs.

**Corollary 1.4 (Strong law of large numbers).** *As $n \to \infty$, $X_n/\log^2 n \to 1/(2m_1)$ almost surely, with $m_1$ defined as in Theorem 1.1.*

Our last main result is a central limit theorem for $\{X_n : n \in \mathbb{N}\}$.

**Theorem 1.5 (Central limit theorem).** *As $n \to \infty$, the sequence*

$$\frac{X_n - (1/(2m_1))\log^2 n}{\sqrt{(m_2/(3m_1^3))\log^3 n}}$$

*weakly converges to the standard normal law, where $m_1$ and $m_2$ are defined as in Theorem 1.1.*

***Remark 1.6.*** The proof of Theorem 1.5 presented in Section 3 draws heavily from coalescent theory and results on random exchangeable partitions. We leave open the question of whether it is possible to deduce the asymptotic normality of $X_n$ from the recursion (2) alone, that is, without using pathwise results available in the coalescent setting.

## 2. Proofs of Theorem 1.1 and Corollary 1.4

**Proof of Theorem 1.1.** For $k \in \mathbb{N}$, set $a_n^{(k)} := \mathbb{E} X_n^k$. By induction on $k$, we will prove the asymptotic expansion

$$a_n^{(k)} = \alpha^k \log^{2k} n + r_k \log^{2k-1} n + \mathrm{O}(\log^{2k-2} n), \qquad k \in \mathbb{N}, \tag{8}$$

where $\alpha := (2m_1)^{-1}$ and

$$r_k := \tfrac{2}{3} k \alpha^{k+1}((2k+1)m_2 + 6cm_1). \tag{9}$$

Recall that $m_1 = m_1^{(b)} = \zeta(2,b)$, $m_2 = m_2^{(b)} = 2\zeta(3,b)$ (see (7)) and $c := -\Psi(b) - 1$.

For $k = 1$, write $a_n$ instead of $a_n^{(1)}$, for simplicity. In view of (2), we have

$$a_1 = 0, \qquad a_n = 1 + \sum_{i=1}^{n-1} a_{n-i} \mathbb{P}\{I_n = i\}, \qquad n \in \{2, 3, \ldots\}. \tag{10}$$

Put $b_n := a_n - \alpha \log^2 n$, $n \in \mathbb{N}$. From (10), it follows that $b_1 = 0$ and

$$\begin{aligned}
b_n &= 1 + \alpha \sum_{i=1}^{n-1}(\log^2(n-i) - \log^2 n)\mathbb{P}\{I_n = i\} + \sum_{i=1}^{n-1} b_{n-i}\mathbb{P}\{I_n = i\} \\
&=: c_n + \sum_{i=1}^{n-1} b_{n-i}\mathbb{P}\{I_n = i\}, \qquad n \in \{2, 3, \ldots\}.
\end{aligned} \tag{11}$$

Using Lemma A.1 (with $k = 1$ and $k = 2$), we get

$$\begin{aligned}
c_n &= 1 + \alpha \sum_{i=1}^{n-1}(\log^2(1 - i/n) + 2\log n \log(1 - i/n))\mathbb{P}\{I_n = i\} \\
&= 1 + \frac{\alpha}{H(n,b)}\left(m_2 + \mathrm{O}\left(\frac{\log^2 n}{n^{b \wedge 1}}\right) + 2\log n\left(-m_1 + \mathrm{O}\left(\frac{\log n}{n^{b \wedge 1}}\right)\right)\right) \\
&= 1 - \frac{\log n}{H(n,b)} + \frac{m_2}{2m_1 H(n,b)} + \mathrm{O}\left(\frac{\log n}{n^{b \wedge 1}}\right)
\end{aligned}$$

and, by (4),

$$c_n = 1 - \frac{H(n,b) + \Psi(b) + 1 + \mathrm{O}(1/n)}{H(n,b)} + \frac{m_2}{2m_1 H(n,b)} + \mathrm{O}\left(\frac{\log n}{n^{b \wedge 1}}\right)$$

$$= \frac{m_2/(2m_1) - \Psi(b) - 1}{H(n,b)} + \mathrm{O}\left(\frac{\log n}{n^{b \wedge 1}}\right) =: \frac{C_1}{H(n,b)} + \mathrm{O}\left(\frac{\log n}{n^{b \wedge 1}}\right).$$

Substituting this relation into (11) yields

$$b_n = \frac{C_1}{H(n,b)} + \mathrm{O}\left(\frac{\log n}{n^{b \wedge 1}}\right) + \sum_{i=1}^{n-1} b_{n-i} \mathbb{P}\{I_n = i\}.$$

Set $d_n := b_n - (C_1/m_1) \log n$, $n \in \mathbb{N}$. Then, $d_1 = 0$ and

$$d_n = \frac{C_1}{H(n,b)} + \frac{C_1}{m_1} \sum_{i=1}^{n-1} \log(1 - i/n) \mathbb{P}\{I_n = i\}$$

$$+ \mathrm{O}\left(\frac{\log n}{n^{b \wedge 1}}\right) + \sum_{i=1}^{n-1} d_{n-i} \mathbb{P}\{I_n = i\}, \qquad n \in \{2, 3, \ldots\}.$$

Another application of Lemma A.1 leads to

$$d_n = \frac{C_1}{H(n,b)} + \frac{C_1}{m_1 H(n,b)} \left(-m_1 + \mathrm{O}\left(\frac{\log n}{n^{b \wedge 1}}\right)\right)$$

$$+ \mathrm{O}\left(\frac{\log n}{n^{b \wedge 1}}\right) + \sum_{i=1}^{n-1} d_{n-i} \mathbb{P}\{I_n = i\}$$

$$= \mathrm{O}\left(\frac{\log n}{n^{b \wedge 1}}\right) + \sum_{i=1}^{n-1} d_{n-i} \mathbb{P}\{I_n = i\}.$$

By Lemma A.2, $d_n = \mathrm{O}(1)$. Therefore, $a_n = \alpha \log^2 n + r_1 \log n + \mathrm{O}(1)$, and we have already proven (8) for $k = 1$.

The induction step from $\{1, \ldots, k\}$ to $k+1$ works as follows. Using (2) and dropping terms of lower orders (which is possible due to the assumption of induction), we get $a_1^{(k+1)} = 0$ and

$$a_n^{(k+1)} = (k+1)\alpha^k \log^{2k} n + (k+1) r_k \log^{2k-1} n$$

$$+ \mathrm{O}(\log^{2k-2} n) + \sum_{j=1}^{n-1} a_{n-j}^{(k+1)} \mathbb{P}\{I_n = j\}, \qquad n \in \{2, 3, \ldots\}.$$

Put $b_n^{(k+1)} := a_n^{(k+1)} - \alpha^{k+1} \log^{2k+2} n$, $n \in \mathbb{N}$. We then have $b_1^{(k+1)} = 0$ and

$$b_n^{(k+1)} = c_n^{(k+1)} + \sum_{j=1}^{n-1} b_{n-j}^{(k+1)} \mathbb{P}\{I_n = j\}, \qquad n \in \{2, 3, \ldots\}, \tag{12}$$

where

$$c_n^{(k+1)} := \alpha^{k+1} \sum_{j=1}^{n-1} (\log^{2k+2}(n-j) - \log^{2k+2} n) \mathbb{P}\{I_n = j\}$$
$$+ (k+1)\alpha^k \log^{2k} n + (k+1)r_k \log^{2k-1} n + O(\log^{2k-2} n).$$

Binomial expansion of $\log^{2k+2}(n-j) = (\log(1-j/n) + \log n)^{2k+2}$ leads to

$$c_n^{(k+1)} = (k+1)\alpha^k \log^{2k} n + (k+1)r_k \log^{2k-1} n + O(\log^{2k-2} n)$$
$$+ \alpha^{k+1} \sum_{j=1}^{n-1} \mathbb{P}\{I_n = j\} \sum_{i=0}^{2k+1} \binom{2k+2}{i} \log^{2k+2-i}(1-j/n) \log^i n$$
$$= (k+1)\alpha^k \log^{2k} n + (k+1)r_k \log^{2k-1} n + O(\log^{2k-2} n)$$
$$+ \alpha^{k+1} \sum_{i=0}^{2k+1} \binom{2k+2}{i} \log^i n \sum_{j=1}^{n-1} \mathbb{P}\{I_n = j\} \log^{2k+2-i}(1-j/n).$$

Using Lemma A.1 gives

$$c_n^{(k+1)} = (k+1)\alpha^k \log^{2k} n + (k+1)r_k \log^{2k-1} n + O(\log^{2k-2} n)$$
$$+ \frac{\alpha^{k+1}}{H(n,b)} \sum_{i=0}^{2k+1} \binom{2k+2}{i} \log^i n \left((-1)^i m_{2k+2-i}^{(b)} + O\left(\frac{\log^{2k+2-i} n}{n^{b \wedge 1}}\right)\right)$$
$$= (k+1)\alpha^k \log^{2k} n + (k+1)r_k \log^{2k-1} n + O(\log^{2k-2} n)$$
$$+ \frac{\alpha^{k+1}}{H(n,b)} \left(-m_1 \binom{2k+2}{2k+1} \log^{2k+1} n + m_2 \binom{2k+2}{2k} \log^{2k} n\right)$$
$$= (k+1)\alpha^k \log^{2k} n \left(1 - \frac{\log n}{H(n,b)}\right)$$
$$+ \left((k+1)r_k + \alpha^{k+1}(2k+1)(k+1)m_2 \frac{\log n}{H(n,b)}\right) \log^{2k-1} n + O(\log^{2k-2} n)$$
$$= (k+1)(r_k + (2k+1)\alpha^{k+1} m_2 - (\Psi(b)+1)\alpha^k) \log^{2k-1} n + O(\log^{2k-2} n)$$
$$=: c_k \log^{2k-1} n + O(\log^{2k-2} n).$$

Plugging the last expression into (12) gives $b_1^{(k+1)} = 0$ and

$$b_n^{(k+1)} = c_k \log^{2k-1} n + \mathrm{O}(\log^{2k-2} n) + \sum_{j=1}^{n-1} b_{n-j}^{(k+1)} \mathbb{P}\{I_n = j\}, \qquad n \in \{2, 3, \ldots\}.$$

Set $e_n^{(k+1)} := b_n^{(k+1)} - C_k \log^{2k+1} n$, $n \in \mathbb{N}$, where $C_k := c_k/((2k+1)m_1)$. The sequence thus defined is given by the recursion

$$e_n^{(k+1)} = c_k \log^{2k-1} n + \mathrm{O}(\log^{2k-2} n)$$

$$+ C_k \sum_{i=1}^{n-1} (\log^{2k+1}(n-i) - \log^{2k+1} n) \mathbb{P}\{I_n = i\}$$

$$+ \sum_{j=1}^{n-1} e_{n-j}^{(k+1)} \mathbb{P}\{I_n = j\}$$

$$= c_k \log^{2k-1} n + \mathrm{O}(\log^{2k-2} n)$$

$$+ C_k \sum_{i=1}^{n-1} \mathbb{P}\{I_n = i\} \sum_{j=0}^{2k} \binom{2k+1}{j} \log^j n \log^{2k+1-j}(1 - i/n)$$

$$+ \sum_{j=1}^{n-1} e_{n-j}^{(k+1)} \mathbb{P}\{I_n = j\}.$$

Again using Lemma A.1 yields

$$e_n^{(k+1)} = c_k \log^{2k-1} n + \mathrm{O}(\log^{2k-2} n)$$

$$+ C_k \frac{\log^{2k} n}{H(n, b)} (2k+1) \left( -m_1 + \mathrm{O}\left(\frac{\log n}{n^{b \wedge 1}}\right) \right) + \sum_{j=1}^{n-1} e_{n-j}^{(k+1)} \mathbb{P}\{I_n = j\}$$

$$= \mathrm{O}(\log^{2k-2} n) + \sum_{j=1}^{n-1} e_{n-j}^{(k+1)} \mathbb{P}\{I_n = j\},$$

by the choice of $C_k$. For sufficiently large $n$, we can choose $M_k > 0$ such that the $\mathrm{O}(\log^{2k-2} n)$ term is dominated by

$$M_k(k\alpha^{k-1} \log^{2k-2} n + kr_{k-1} \log^{2k-3} n + \mathrm{O}(\log^{2k-4} n)).$$

Therefore, for large $n$, $e_n^{(k+1)} \le M_k a_n^{(k)}$. By the assumption of induction, $a_n^{(k)} = \mathrm{O}(\log^{2k} n)$. Therefore, $e_n^{(k+1)} = \mathrm{O}(\log^{2k} n)$. Setting $r_{k+1} := C_k = c_k/((2k+1)m_1)$, we obtain

$$a_n^{(k+1)} = \alpha^{k+1} \log^{2k+2} n + r_{k+1} \log^{2k+1} n + \mathrm{O}(\log^{2k} n).$$

The sequence $\{r_k : k \in \mathbb{N}\}$ satisfies the recursion

$$r_{k+1} = \frac{k+1}{(2k+1)m_1}(r_k + (2k+1)\alpha^{k+1}m_2 - (\Psi(b)+1)\alpha^k)$$

with initial condition

$$r_1 = \frac{m_2/(2m_1) - \Psi(b) - 1}{m_1} = \frac{\zeta(3,b)/\zeta(2,b) - \Psi(b) - 1}{\zeta(2,b)}.$$

The unique solution of this recursion is given by (9). The proof of Theorem 1.1 is thus complete. □

**Proof of Corollary 1.4.** For $n \in \mathbb{N}$ and $\varepsilon > 0$, set $A_n(\varepsilon) := \{|X_n - \mathbb{E}X_n| \geq \varepsilon \mathbb{E}X_n\}$. By Chebyshev's inequality, $\mathbb{P}\{A_n(\varepsilon)\} \leq \mathbb{D}X_n/(\varepsilon \mathbb{E}X_n)^2$. From Theorem 1.1, it follows that

$$\frac{\mathbb{D}X_n}{(\mathbb{E}X_n)^2} = \frac{4m_2}{3m_1}\frac{1}{\log n} + O\left(\frac{1}{\log^2 n}\right).$$

Therefore, replacing $n$ by $n_k := \lfloor \exp(k^2) \rfloor$, it follows that $\sum_{k=1}^{\infty} \mathbb{P}\{A_{n_k}(\varepsilon)\} < \infty$ and hence $X_{n_k}/\mathbb{E}X_{n_k} \to 1$ almost surely as $k \to \infty$, by the Borel–Cantelli lemma. Thus, we have already verified the result along the subsequence $\{n_k : k \in \mathbb{N}\}$. For each integer $n \geq n_1$, there exists a unique index $k = k(n) \in \mathbb{N}$ such that $n_k \leq n < n_{k+1}$. By its definition, the sequence $\{X_n : n \in \mathbb{N}\}$ is almost surely non-decreasing. Therefore, the corollary follows from the standard sandwich argument

$$\frac{X_{n_k}}{\mathbb{E}X_{n_k}}\frac{\mathbb{E}X_{n_k}}{\mathbb{E}X_{n_{k+1}}} \leq \frac{X_n}{\mathbb{E}X_n} \leq \frac{X_{n_{k+1}}}{\mathbb{E}X_{n_{k+1}}}\frac{\mathbb{E}X_{n_{k+1}}}{\mathbb{E}X_{n_k}} \qquad \text{almost surely}$$

and from $\mathbb{E}X_{n_k}/\mathbb{E}X_{n_{k+1}} \sim \log^2 n_k / \log^2 n_{k+1} \sim k^4/(k+1)^4 \to 1$. □

## 3. Proof of Theorem 1.5

We will use Theorem 2.1 of Neininger and Rüschendorf [13], which is written down below in a modified form suggested by Gnedin, Pitman and Yor [7], Theorem 10. In the following, for random variables $X$, we use the notation $\|X\|_3 := (\mathbb{E}(|X|^3))^{1/3}$.

**Proposition 3.1.** *Assume that a random sequence $\{U_n : n \in \mathbb{N}\}$ of scalar random variables satisfies the recursion*

$$U_n \stackrel{d}{=} U_{J_n} + V_n, \qquad n \in \{n_0, n_0 + 1, \ldots\}, \tag{13}$$

*for some $n_0 \in \mathbb{N}$, where $(J_n, V_n)$ is independent of $\{U_n : n \geq n_0\}$, $J_n$ takes values in $\{0, 1, \ldots, n\}$ and $\mathbb{P}\{J_n = n\} < 1$ for each integer $n \geq n_0$. Suppose, further, that $\|U_n\|_3 < \infty$ and that for some constants $C > 0$ and $\alpha > 0$, the following three conditions hold:*

(i) $\limsup_{n\to\infty} \mathbb{E}\log(\frac{J_n\vee 1}{n}) < 0$ *and* $\sup_{n\geq 2} \|\log(\frac{J_n\vee 1}{n})\|_3 < \infty$;
(ii) *for some* $\lambda \in [0, 2\alpha)$ *and some* $\kappa > 0$, *as* $n \to \infty$,

$$\|V_n - \mu_n + \mu_{J_n}\|_3 = \mathrm{O}(\log^\kappa n), \qquad \mathbb{D}U_n = C\log^{2\alpha} n + \mathrm{O}(\log^\lambda n),$$

*where* $\mu_n := \mathbb{E}U_n$;
(iii) $\alpha > 1/3 + \max(\kappa, \lambda/2)$.

*Then, as* $n \to \infty$, *the sequence* $(U_n - \mu_n)/(\sqrt{C}\log^\alpha n)$ *weakly converges to the standard normal law.*

**Remark 3.2.** The recursion (2) is of the form (13) with random indices $J_n := n - I_n$, where $I_n$ has distribution (3). By Lemma A.1 and (4),

$$\mathbb{E}\log\left(\frac{J_n}{n}\right) = \sum_{i=1}^{n-1} \log\left(1 - \frac{i}{n}\right)\mathbb{P}\{I_n = i\} \sim -\frac{m_1^{(b)}}{\log n}.$$

Therefore, $\lim_{n\to\infty} \mathbb{E}\log(J_n/n) = 0$. In particular, the first part of condition (i) in Proposition 3.1 is not satisfied. Hence, Proposition 3.1 is not applicable to the recursion (2).

Fix any $T > 0$. The total number $X_n$ of collisions is the sum of the numbers of collisions occurring during the time intervals $[0, T)$ (denote this by $X_n(T)$) and $[T, \infty)$ (denote this by $\widehat{X}_n(T)$). Since the coalescent is a Markov process, $\widehat{X}_n(T) \stackrel{d}{=} X'_{|\varrho_n\Pi_T|}$, where $(J_n, V_n) := (|\varrho_n\Pi_T|, X_n(T))$ is independent of $\{X'_n : n \in \mathbb{N}\}$, an independent copy of $\{X_n : n \in \mathbb{N}\}$. Thus, we have proven that $\{X_n : n \in \mathbb{N}\}$ satisfies another recursion of the form (13), namely

$$X_n \stackrel{d}{=} X_{|\varrho_n\Pi_T|} + X_n(T). \tag{14}$$

**Proof of Theorem 1.5.** Let us prove that the recursion (14) satisfies all of the conditions of Proposition 3.1.

Since $X_n \leq n - 1$ almost surely, $\|X_n\|_3 < \infty$.

Recall that $X_n(T)$ is the number of jumps of the process $\{\varrho_n\Pi_t : t \in [0, T)\}$. If $\Lambda$ has no atom at the origin, then any $\Lambda$-coalescent can be constructed from a Poisson point process (see page 1872 in [15]). From this construction, it follows that with probability one, $X_n(T)$ is bounded from above by a random variable with Poisson distribution with parameter $Tg_n$. By (5), $Tg_n \sim (T/B(2,b))\log n$. Therefore,

$$\|X_n(T)\|_3 = \mathrm{O}(\log n), \qquad n \to \infty. \tag{15}$$

Let $Q(T) := \{\widehat{f}_i(T) : i \in \mathbb{N}\}$ be the decreasing rearrangement of the asymptotic frequencies of the random exchangeable partition $\Pi_T$. According to Remark 1.3, $1 - \sum_{i=1}^\infty \widehat{f}_i(T) = \mathrm{e}^{-\widehat{S}_T}$. The elements of the set $Q(T) \cup \{1 - \sum_{i=1}^\infty \widehat{f}_i(T)\}$ are the lengths of the intervals (from left to right) comprising the partition of $[0, 1]$. Let $U_1, \ldots, U_n$ be independent random variables (points), uniformly distributed on $[0, 1]$ and independent of the lengths of the intervals. Let $W_{n,i}(T)$ be the number of points falling in the interval of length $\widehat{f}_i(T)$. Set

$$\eta_n(T) := |\{i \in \{1, \ldots, n\} : U_i > 1 - e^{-\widehat{S}_T}\}|,$$

$$\zeta_n(T) := |\{i \geq 1 : W_{n,i}(T) > 0\}|, \qquad \theta_n(T) := \zeta_n(T) + 1_{\{\eta_n(T) > 0\}}.$$

From the paintbox construction [12] of a random exchangeable partition, it follows that

$$|\varrho_n \Pi_T| \stackrel{d}{=} \zeta_n(T) + \eta_n(T).$$

Arguing in the same way as on page 592 in [7], we conclude that as $n \to \infty$, $\eta_n(T)/n \to e^{-\widehat{S}_T}$ almost surely, which easily implies that

$$\lim_{n \to \infty} \left( -\log\left(\frac{\eta_n(T) \vee 1}{n}\right) \right) = \widehat{S}_T \tag{16}$$

almost surely and that for each $k \in \mathbb{N}$,

$$\lim_{n \to \infty} \mathbb{E}\left| \left( \log\left(\frac{\eta_n(T) \vee 1}{n}\right) \right)^k \right| = \mathbb{E}\widehat{S}_T^k. \tag{17}$$

Note that, in view of (7), the right-hand side is finite for each $k \in \mathbb{N}$. Interpreting the intervals as "boxes" and the points as "balls", the $\theta_n(T)$ is just the number of occupied boxes in the classical multinomial occupancy scheme. From the results on page 152 in [5], it follows that $\lim_{n \to \infty} n^{-1} \mathbb{E}(\theta_n(T) | \widehat{f}_i(T) : i \in \mathbb{N}) = 0$ almost surely. This fact, together with Proposition 2 of the same reference (see also Theorem 8 in [11]), leads to $\lim_{n \to \infty} \theta_n(T)/n = 0$ almost surely conditionally on $\{\widehat{f}_i(T) : i \in \mathbb{N}\}$ and, hence, unconditionally. The latter implies that $\lim_{n \to \infty} |\varrho_n \Pi_T|/n = e^{-\widehat{S}_T}$ almost surely and, hence,

$$\lim_{n \to \infty} \left( -\log\left(\frac{|\varrho_n \Pi_T|}{n}\right) \right) = \widehat{S}_T \tag{18}$$

almost surely. Since

$$-\log\left(\frac{|\varrho_n \Pi_T|}{n}\right) \leq -\log\left(\frac{\eta_n(T) \vee 1}{n}\right)$$

almost surely, (16)–(18) together imply that for each $k \in \mathbb{N}$,

$$\lim_{n \to \infty} \mathbb{E}\left| \left( \log\left(\frac{|\varrho_n \Pi_T|}{n}\right) \right)^k \right| = \mathbb{E}\widehat{S}_T^k, \tag{19}$$

by a variant of Fatou's lemma sometimes called Pratt's lemma (see [16]).

Condition (i) of Proposition 3.1 follows from (19). The estimate $\|\mu_n - \mu_{J_n}\|_3 = \mathrm{O}(\log n)$ follows from Theorem 1.1 and (19). In view of this observation, (15) and Theorem 1.1, (ii) holds with $\kappa = 1$, $\alpha = 3/2$ and $\lambda = 2$. Therefore, (iii) also holds. □

# Appendix

The proof of Theorem 1.1 relies on the two following technical results.

**Lemma A.1.** *For all $k \in \mathbb{N}$ and $b > 0$, as $n \to \infty$,*

$$\left| H(n,b) \sum_{i=1}^{n-1} \mathbb{P}\{I_n = i\} \left(-\log\left(1 - \frac{i}{n}\right)\right)^k - m_k^{(b)} \right| = \mathrm{O}\left(\frac{\log^k n}{n^{b \wedge 1}}\right), \qquad (20)$$

*where $H(n,b)$ is the function defined after (3) and $m_k^{(b)} = k!\zeta(k+1,b)$ is the $k$th moment (see (7)) of the Lévy measure (6).*

**Proof.** We first prove that

$$J_n(b,k) := \left| \sum_{i=1}^{n-1} \left(1 - \frac{i}{n}\right)^{b-1} \frac{1}{i} \left(-\log\left(1 - \frac{i}{n}\right)\right)^k - m_k^{(b)} \right| = \mathrm{O}\left(\frac{\log^k n}{n^{b \wedge 1}}\right) \qquad (21)$$

and that

$$L_n(b,k) := \left| \sum_{i=1}^{n-1} \left(1 - \frac{i}{n}\right)^{b-1} \frac{1}{i+1} \left(-\log\left(1 - \frac{i}{n}\right)\right)^k - m_k^{(b)} \right| = \mathrm{O}\left(\frac{\log^k n}{n^{b \wedge 1}}\right). \qquad (22)$$

Fix $k \in \mathbb{N}$. For $b > 1$, introduce the continuous non-negative function $f_b : [0,1] \to \mathbb{R}$ via $f_b(x) := x^{-1}(1-x)^{b-1}(-\log(1-x))^k$ for $x \in (0,1)$, $f_b(0) := 1_{\{k=1\}}$ and $f_b(1) := 0$. Pick some $\delta \in (0,1)$ such that $f_b$ is non-increasing on $[\delta, 1]$. Then,

$$\left| \frac{1}{n} \sum_{i=[n\delta]+1}^{n-1} f_b\left(\frac{i}{n}\right) - \int_\delta^1 f_b(x)\, \mathrm{d}x \right|$$

$$= \left| \sum_{i=[n\delta]+1}^{n-1} \int_{i/n}^{(i+1)/n} \left( f_b\left(\frac{i}{n}\right) - f_b(x) \right) \mathrm{d}x - \int_\delta^{([n\delta]+1)/n} f_b(x)\, \mathrm{d}x \right|$$

$$\leq \sum_{i=[n\delta]+1}^{n-1} \int_{i/n}^{(i+1)/n} \left( f_b\left(\frac{i}{n}\right) - f_b\left(\frac{i+1}{n}\right) \right) \mathrm{d}x + \int_\delta^{([n\delta]+1)/n} f_b(x)\, \mathrm{d}x$$

$$= \mathrm{O}\left(\frac{1}{n}\right).$$

It is easily checked that $f_b$ is continuously differentiable on $(0, \delta)$ with $\sup_{0 < x < \delta} |f_b'(x)| < \infty$. Therefore, exploiting the mean value theorem for differentiable functions, we have

$$\left| \frac{1}{n} \sum_{i=1}^{[n\delta]} f_b\left(\frac{i}{n}\right) - \int_0^\delta f_b(x)\, \mathrm{d}x \right| = \mathrm{O}\left(\frac{1}{n}\right).$$

Combining these two pieces and using the equality $m_k^{(b)} = \int_0^1 f_b(x)\,\mathrm{d}x$, we get $J_n(b,k) = \mathrm{O}(1/n)$, which is more than we need.

Assuming that $b \in (0,1]$, an application of the previous result to the function $f_{b+1}$, which satisfies

$$f_{b+1}(x) = \frac{(1-x)^{b-1}(-\log(1-x))^k}{x} - (1-x)^{b-1}(-\log(1-x))^k$$

for $x \in (0,1)$, yields

$$\left| \sum_{i=1}^{n-1} \left( \frac{(1-i/n)^{b-1}(-\log(1-i/n))^k}{i} - \frac{(1-i/n)^{b-1}(-\log(1-i/n))^k}{n} \right) \right.$$

$$\left. - \int_0^1 f_{b+1}(x)\,\mathrm{d}x \right| = \mathrm{O}\!\left(\frac{1}{n}\right). \tag{23}$$

Note that $\int_0^1 f_{b+1}(x)\,\mathrm{d}x = m_k^{(b)} - k!/b^{k+1}$.

For all $n \in \mathbb{N}$ with $b \log n \geq 1$, we now use the inequalities

$$\frac{1}{n}\sum_{i=1}^{n-1}\left(\frac{i}{n}\right)^{b-1}\left(-\log\left(\frac{i}{n}\right)\right)^k \geq \int_{1/n}^1 x^{b-1}(-\log x)^k\,\mathrm{d}x = \frac{k!}{b^{k+1}}\left(1 - n^{-b}\sum_{i=0}^{k}\frac{(b\log n)^i}{i!}\right)$$

$$\geq \frac{k!}{b^{k+1}} - k!\frac{\log^k n}{bn^b}\sum_{i=0}^{k}\frac{1}{i!} \geq \frac{k!}{b^{k+1}} - k!\mathrm{e}\frac{\log^k n}{bn^b}$$

to conclude that, as $n \to \infty$,

$$\left| \frac{1}{n}\sum_{i=1}^{n-1}\left(1 - \frac{i}{n}\right)^{b-1}\left(-\log\left(1 - \frac{i}{n}\right)\right)^k - \frac{k!}{b^{k+1}} \right| = \mathrm{O}\!\left(\frac{\log^k n}{n^b}\right).$$

Combining this estimate with (23) yields (21).

Let us now prove (22). If $k \in \mathbb{N} \setminus \{1\}$, then

$$0 \leq M_n(b,k)$$

$$:= \sum_{i=1}^{n-1}\frac{(1-i/n)^{b-1}(-\log(1-i/n))^k}{i} - \sum_{i=1}^{n-1}\frac{(1-i/n)^{b-1}(-\log(1-i/n))^k}{i+1}$$

$$\leq \sum_{i=1}^{n-1}\frac{(1-i/n)^{b-1}(-\log(1-i/n))^k}{i^2}$$

$$\sim \frac{1}{n}\int_0^1 \frac{(1-x)^{b-1}(-\log(1-x))^k}{x^2}\,\mathrm{d}x$$

and the last integral is finite. Therefore, $M_n(b,k) = O(1/n)$, which, together with (21), proves (22) under the current assumptions.

If $k = 1$, then

$$\begin{aligned}
0 &\leq M_n(b,k) \\
&\leq n^{(1-b)\vee 0} \sum_{i=1}^{n-1} \frac{-\log(1-i/n)}{i^2} = n^{(1-b)\vee 0} \sum_{i=1}^{n-1} \frac{1}{i^2} \sum_{j=1}^{\infty} \frac{(i/n)^j}{j} \\
&\leq n^{(1-b)\vee 0} \sum_{i=1}^{n-1} \frac{1}{i^2} \sum_{j=1}^{\infty} \left(\frac{i}{n}\right)^j = n^{(1-b)\vee 0} \sum_{i=1}^{n-1} \frac{1}{i^2} \frac{i/n}{1-i/n} \\
&= n^{(1-b)\vee 0} \sum_{i=1}^{n-1} \frac{1}{i(n-i)} = n^{(1-b)\vee 0} \frac{1}{n} \sum_{i=1}^{n-1} \left(\frac{1}{i} + \frac{1}{n-i}\right) \sim \frac{2\log n}{n^{b\wedge 1}}.
\end{aligned}$$

This relation, together with (21), proves (22).

For $b=1$, the left-hand side of (22) coincides with that of (20). Thus, we only have to check (20) for $b \neq 1$. To this end, keeping in mind (21) and (22), it suffices to show that

$$\left| \sum_{i=1}^{n-1} \left( \frac{\Gamma(n-i+b-1)\Gamma(n+1)}{\Gamma(n-i)\Gamma(n+b)} - \left(1 - \frac{i}{n}\right)^{b-1} \right) \frac{1}{i+1} \left( -\log\left(1 - \frac{i}{n}\right) \right)^k \right| \\
= O\left( \frac{\log^k n}{n^{b\wedge 1}} \right). \tag{24}$$

First, we will prove that for any $b > 0$, there exists a constant $M > 0$ such that for all $n \in \mathbb{N}$ and all $j \in \{1, \ldots, n-1\}$,

$$\left| \frac{\Gamma(n-j+b-1)\Gamma(n+1)}{\Gamma(n-j)\Gamma(n+b)} - \left(1 - \frac{j}{n}\right)^{b-1} \right| \leq \frac{M}{n} \left(1 - \frac{j}{n}\right)^{b-2} \tag{25}$$

or, equivalently,

$$\left| \frac{\Gamma(j+b-1)\Gamma(n+1)}{\Gamma(j)\Gamma(n+b)} - \left(\frac{j}{n}\right)^{b-1} \right| \leq \frac{M}{n} \left(\frac{j}{n}\right)^{b-2}. \tag{26}$$

The subsequent argument relies on the following inequality (see (6.1.47) in [1]). For $c, d > -1$, there exists $M_{c,d} > 0$ such that for all $n \in \mathbb{N}$,

$$\left| \frac{\Gamma(n+c)}{\Gamma(n+d)} - n^{c-d} \right| \leq M_{c,d} n^{c-d-1}.$$

(26) now follows from the chain of inequalities

$$\left| \frac{\Gamma(j+b-1)\Gamma(n+1)}{\Gamma(j)\Gamma(n+b)} - \left(\frac{j}{n}\right)^{b-1} \right|$$

$$= \left|\left(\frac{\Gamma(j+b-1)}{\Gamma(j)} - j^{b-1}\right)\frac{\Gamma(n+1)}{\Gamma(n+b)} + \frac{\Gamma(n+1)}{\Gamma(n+b)}j^{b-1} - \left(\frac{j}{n}\right)^{b-1}\right|$$

$$\leq \frac{\Gamma(n+1)}{\Gamma(n+b)}\left|\frac{\Gamma(j+b-1)}{\Gamma(j)} - j^{b-1}\right| + j^{b-1}\left|\frac{\Gamma(n+1)}{\Gamma(n+b)} - n^{1-b}\right|$$

$$\leq \frac{\Gamma(n+1)}{\Gamma(n+b)}M_{b-1,0}j^{b-2} + j^{b-1}M_{1,b}n^{-b}$$

$$\leq \left|\frac{\Gamma(n+1)}{\Gamma(n+b)} - n^{1-b}\right|M_{b-1,0}j^{b-2} + n^{1-b}M_{b-1,0}j^{b-2} + j^{b-1}M_{1,b}n^{-b}$$

$$\leq \frac{M_{1,b}M_{b-1,0}}{n^2}\left(\frac{j}{n}\right)^{b-2} + \frac{M_{b-1,0}}{n}\left(\frac{j}{n}\right)^{b-2} + \frac{M_{1,b}}{n}\left(\frac{j}{n}\right)^{b-1}$$

$$\leq \frac{M}{n}\left(\frac{j}{n}\right)^{b-2},$$

where $M := M_{1,b}M_{b-1,0} + M_{b-1,0} + M_{1,b}$. Plugging (25) into the left-hand side of (24) gives

$$\left|\sum_{i=1}^{n-1}\left(\frac{\Gamma(n-i+b-1)\Gamma(n+1)}{\Gamma(n-i)\Gamma(n+b)} - \left(1-\frac{i}{n}\right)^{b-1}\right)\frac{1}{i+1}\left(-\log\left(1-\frac{i}{n}\right)\right)^k\right|$$

$$\leq \frac{M}{n}\sum_{i=1}^{n-1}\left(1-\frac{i}{n}\right)^{b-2}\frac{1}{i+1}\left(-\log\left(1-\frac{i}{n}\right)\right)^k =: Q_n(b,k).$$

For $b > 1$, the function $x \mapsto x^{-1}(1-x)^{b-2}\log^k(1-x)$ is integrable on $[0,1]$, which implies that the latter sum is bounded and the right-hand side in (24) is $\mathrm{O}(1/n)$. If $b \in (0,1)$, then noting that the function $x \mapsto x^{-1}(-\log(1-x))^k$ is non-decreasing on $(0,1)$, we conclude that for $n \in \{2, 3, \ldots\}$,

$$Q_n(b,k) = \frac{M}{n^b}\sum_{i=1}^{n-1}(n-i)^{b-2}\frac{1}{(i+1)/n}\left(-\log\left(1-\frac{i}{n}\right)\right)^k$$

$$\leq \frac{M}{n^b}\sum_{i=1}^{n-1}(n-i)^{b-2}\frac{1}{i/n}\left(-\log\left(1-\frac{i}{n}\right)\right)^k$$

$$\leq \frac{2M\log^k n}{n^b}\sum_{i=1}^{n-1}(n-i)^{b-2} = \mathrm{O}\left(\frac{\log^k n}{n^b}\right).$$

Thus, (24) is established and the proof is complete. □

**Lemma A.2.** *Fix $k \in \mathbb{N}$ and $b > 0$, and suppose that $\{a_n : n \in \mathbb{N}\}$ is some sequence satisfying $a_n = \mathrm{O}(n^{-b} \log^k n)$. If the sequence $\{v_n : n \in \mathbb{N}\}$ is defined recursively by*

$$v_1 := 0, \qquad v_n := a_n + \sum_{i=1}^{n-1} v_{n-i} \mathbb{P}\{I_n = i\}, \qquad n \in \{2, 3, \ldots\},$$

*where $\mathbb{P}\{I_n = k\}$ is defined in (3), then $v_n = \mathrm{O}(1)$.*

**Proof.** Since $\mathbb{E} I_n \sim n/(b \log n)$, there exists an $M > 0$ such that for all $n \in \{2, 3, \ldots\}$,

$$\frac{b}{2n^{1+b/2}} \mathbb{E} I_n \geq \frac{M \log^k n}{n^b}. \qquad (27)$$

It suffices to prove the following. If

$$u_1 := 0, \qquad u_n = \frac{M \log^k n}{n^b} + \sum_{i=1}^{n-1} u_{n-i} \mathbb{P}\{I_n = i\}, \qquad n \in \{2, 3, \ldots\},$$

with $M$ defined in (27), then

$$u_n \leq 2 - n^{-b/2} \qquad \text{for all } n \in \mathbb{N}. \qquad (28)$$

We will use induction. For $n = 1$, (28) is obviously satisfied as $u_1 = 0$. Assume (28) holds for all $n \in \{1, \ldots, m-1\}$. Then,

$$u_m \leq \frac{M \log^k m}{m^b} + \sum_{i=1}^{m-1} (2 - (m-i)^{-b/2}) \mathbb{P}\{I_m = i\}.$$

We will now verify that the right-hand side of the latter inequality is less than or equal to $2 - m^{-b/2}$ or, equivalently, that

$$\sum_{i=1}^{m-1} ((m-i)^{-b/2} - m^{-b/2}) \mathbb{P}\{I_m = i\} \geq \frac{M \log^k m}{m^b}.$$

The inequality $(1 - x)^{-a} \geq 1 + ax$, $x \in (0, 1)$, $a > 0$ yields

$$\sum_{i=1}^{m-1} ((m-i)^{-b/2} - m^{-b/2}) \mathbb{P}\{I_m = i\}$$

$$= m^{-b/2} \sum_{i=1}^{m-1} ((1 - i/m)^{-b/2} - 1) \mathbb{P}\{I_m = i\}$$

$$\geq \frac{b}{2m^{1+b/2}} \mathbb{E} I_m \geq \frac{M \log^k m}{m^b},$$

by (27). $\qquad \square$

# Acknowledgements

The first author was supported by the German Research Foundation DFG, Project 436UKR 113/93/0-1. The authors would like to thank the referees for their careful reading and for their suggestions leading to a significant improvement of the style of the manuscript.